\documentclass[12pt,reqno,a4paper]{amsart}
\usepackage[utf8]{inputenc}
\usepackage{amsmath,pdfsync,verbatim,graphicx,epstopdf,enumerate,latexsym,amssymb}
\usepackage[margin=1in]{geometry}
\usepackage[usenames]{color}

\usepackage[abbrev,nobysame,alphabetic]{amsrefs}
\usepackage{mathtools}  
\usepackage{thmtools}
\usepackage{thm-restate}
\mathtoolsset{showonlyrefs}

\usepackage[colorlinks=true]{hyperref}
\usepackage{cancel}
\usepackage{accents}
\usepackage[framemethod=tikz]{mdframed}
\hypersetup{linkcolor=blue,citecolor=blue}
\usepackage{amsmath,pdfsync,verbatim,graphicx,epstopdf,enumerate}

\usepackage{mathtools}
\usepackage{amsthm}
\newtheorem{theorem}{Theorem}[section]

\newtheorem{proposition}[theorem]{Proposition}

\newtheorem{remark}[theorem]{Remark}

\newcommand{\R}{\mathbb{R}}

\newcommand{\Lc}{\mathcal{L}}

\newcommand{\Oc}{\mathcal{O}}
\newcommand{\Cb}{\mathbb{C}}

\usepackage{amsmath,pdfsync,verbatim,graphicx,epstopdf,enumerate}

\renewcommand{\O}{\Omega}

\newcommand{\C}{\mathbb{C}}

\newcommand{\p}{\partial}

\newcommand{\PO}{\partial\Omega}
\renewcommand{\O}{\Omega}

\title[Local data inverse problem for convection-diffusion equation]{Inverse boundary value problem for the convection-diffusion equation with local data}

\author[Kumar and Purohit]{Pranav Kumar and Anamika Purohit}
 \address{Department of Mathematics, Indian Institute of Science Education and Research, Bhopal - 462066.
 \newline
 E-mail:{\tt \ pranav19@iiserb.ac.in}}
 \address {Department of Mathematics, Indian Institute of Technology, Gandhinagar, Gujarat - 382355.
 \newline
 E-mail:{\tt \ anamika.purohit@iitgn.ac.in}}

\begin{document}
\begin{abstract}
We study a local data inverse problem for the time-dependent convection-diffusion equation in a bounded domain where a part of the boundary is treated to be inaccessible. Up on assuming the inaccessible part to be flat, we seek for the unique determination of the time-dependent convection and the density terms from the knowledge of the boundary data measured outside the inaccessible part. In the process, we show that there is a natural gauge in the perturbations, and we prove that this is the only obstruction in the uniqueness result.
\end{abstract}

\maketitle
\vspace{-5mm}
\vspace{2mm}
\section{Introduction and the main result} 
Let $\Omega\subset\R^{n}$ be a bounded simply connected open set with smooth boundary $\PO$, where $n\geq3$. For $T>0$, let $Q:= (0,T)\times \Omega$ and $\Sigma:= (0,T)\times \partial\Omega$. Let $\Gamma_{0}$ be an arbitrary non-empty open subset of the boundary $\partial\Omega$ and $\Gamma :=\partial\Omega\setminus \Gamma_{0}$, then we denote $\Sigma^{\sharp}:= (0,T)\times \Gamma$. Given a vector field $A(t,x)\in C^{\infty}(\overline{Q}; \Cb^n)$ and function $q(t,x) \in L^{\infty}(Q; \C)$, we consider the following initial boundary value problem (IBVP)

\begin{equation}\label{IBVP}
    \begin{aligned}
\begin{cases}
\Lc_{A,q} u(t,x)=0 &\text { in } Q,\\
u(0,x)=0 &\text{ in }\O,\\
u(t,x)=f(t,x) &\text{ on } \Sigma,
\end{cases}
    \end{aligned}
\end{equation}
where 
\begin{align}\label{Convection_diffusion_operator}
    \Lc_{A,q}:=&\partial_{t}-\sum_{j=1}^{n}\left(\partial_{j}+{A_{j}(t,x)}\right)^{2}+q(t,x).
\end{align}
Here, $A$ is the convection term, $q$ is the density coefficient with constant diffusion, and the IBVP \eqref{IBVP} is known as a convection-diffusion equation. This equation appears in chemical engineering in describing the movement of macro-particles; see \cite{stocker2011introduction} and in the probabilistic study of the diffusion process (Fokker-Planck and Kolmogorov equations).

We begin with discussing the well-posedness of forward problem \eqref{IBVP}. Following \cite{caro2019determination}, we define the following spaces:
\begin{equation}\label{Domain of Dirichlet data}
\mathcal{K}_0(\Sigma):= \{f|_{\Sigma} : f \in L^2(0,T;H^1(\Omega)) \cap H^1(0,T;H^{-1}(\Omega)) \text{ and } f(0,x)=0, \text{ for } x\in \Omega\}, 
\end{equation}
and 
\begin{equation}\label{Space H-T}
\mathcal{H}_T(\Sigma):= \{g|_{\Sigma} : g \in H^1(0,T;H^1(\Omega))  \text{ and } g(T,x)=0, \text{ for }  x\in \Omega\}.
\end{equation} 
Following \cite{caro2019determination}*{Proposition 2.2}, (see also \cite{lions1968problèmes}*{Theorem 4.1, Chapter 3}), one can show  that for any $f\in \mathcal{K}_{0}(\Sigma)$, there exists a unique solution $u\in L^2(0,T;H^1(\Omega))\cap H^1(0,T;H^{-1}(\Omega))$ to IBVP \eqref{IBVP}. Therefore, we can define the Dirichlet--to--Neumann (DN) map $\Lambda_{A,q}:\mathcal{K}_0(\Sigma)\to \mathcal{H}^*_T(\Sigma)$ as
 \begin{align}\label{D-N map}
    \Lambda_{A,q}(f)=  \left.(\partial_\nu u + 2(\nu\cdot A) u)\right|_{\Sigma}, \quad\quad  f\in \mathcal{K}_{0}(\Sigma),
\end{align}
where $\mathcal{H}^*_T(\Sigma)$ denotes the dual of space $\mathcal{H}_T(\Sigma)$ and $u$ solves \eqref{IBVP} and $\nu$ denotes the unit outward normal vector to the boundary $\partial\Omega$.
We also define the local DN map 
$\Lambda^{\Sigma^{\sharp}}_{A,q}:\mathcal{K}_0(\Sigma)\to \mathcal{H}^*_T(\Sigma)$ by
 \begin{align}\label{local-DN}
    \Lambda_{A,q}^{\Sigma^{\sharp}}(f)=  \left.(\partial_\nu u + 2(\nu\cdot A) u)\right|_{\Sigma^{\sharp}}, \quad\quad  f\in \mathcal{K}_{0}(\Sigma), \quad \mathrm{supp}(f(t,\cdot)|_{\Sigma})\subset \Gamma, \quad \text { for each $t\in (0,T)$.
    }
\end{align}
To the best of the author's knowledge, the unique determination of the convection term $A$ and density coefficient $q$ from the knowledge of the local DN map \eqref{local-DN} have not been studied earlier.
\begin{remark}[Gauge invariance]\label{Gauge-Invariance}
As discussed in \cites{caro2019determination, mishra2022inverse}, from the DN map, one can only recover the convection term $A$ and the density coefficient $q$ up to a gauge invariance.
Let $u$ be the solution of \eqref{IBVP} and $\Psi \in C^{\infty}(\overline{Q})$ such that $\Psi|_{\Sigma} = \partial_{\nu}\Psi|_{\Sigma} = 0$ and $\widetilde{A}=A+\nabla_{x} \Psi$, $\widetilde{q}=q+\partial_{t}\Psi$. Then, $\widetilde{u} :=  e^{-\Psi}u$, satisfies \eqref{IBVP} with $A=\widetilde{A}$ and $q=\widetilde{q}$.
Moreover, the DN map corresponding to both equations remains the same, that is, $\Lambda_{A,q}=\Lambda_{\widetilde{A},\widetilde{q}}$ but $(A, q)\neq (\widetilde{A},\widetilde{q})$ as long as $\Psi \neq 0$.
\end{remark}    
 
The following is the main result of the article, where we prove that the gauge above is the only obstruction in recovering the perturbations from the knowledge of the local DN map.
\begin{theorem}\label{Theorem_hyperplane}
    Let $Q=(0, T)\times \Omega$ where $T>0$ and $\Omega \subset \left\{x\in \mathbb{R}^{n}: x_{n}>0\right\}$, $n\geq 3,$ be a bounded simply connected domain with smooth boundary $\PO$, and let $\Gamma_{0} = \partial\Omega \cap \left\{x\in {\R}^{n}: x_{n}=0\right\}\neq \emptyset$, and $\Gamma:= \partial \Omega \setminus \Gamma_{0} $. Let $A^{(j)}\in C^{\infty}(\overline{Q}; \C^{n})$, $q_{j}\in L^{\infty}(Q)$ such that \eqref{IBVP} is well-posed for $j=1,2$. Let $\Lambda_{A^{(j)}, q_{j}}^{\Sigma^{\sharp}}$ be the local Dirichlet-to-Neumann map corresponding to $u_j$ for $j=1,2$.
    If 
    \begin{align}
    \label{eq:Equal_DN_map}
    \Lambda_{A^{(1)}, q_{1}}^{\Sigma^{\sharp}}(f)=\Lambda_{A^{(2)}, q_{2}}^{\Sigma^{\sharp}}(f), \quad \text {for all } f\in \mathcal{K}_{0}(\Sigma) \quad \text { with } \quad \mathrm{supp}(f(t,\cdot)|_{\Sigma})\subset \Gamma, \text { for each $t\in (0,T)$}.
    \end{align}
    Then there exists a function $\Psi\in C^{\infty}(\overline{Q})$ with $\Psi|_{\Sigma} = \partial_{\nu}\Psi|_{\Sigma} = 0$ such that $A^{(2)}(t,x)-A^{(1)}(t,x) = \nabla_x\Psi(t,x)$ and $q_2(t,x) - q_1(t,x) = \partial_t\Psi(t,x)$ in $Q$ provided $A^{(1)}(t,x)=A^{(2)}(t,x)$ on $\Sigma$.
\end{theorem}

The problem concerning the unique recovery of the coefficients for parabolic partial differential equations (PDEs) from boundary measurements have attracted considerable attention recently. Motivated by the seminal work \cite{Global_uniqueness}, Isakov in \cite{isakov1991completeness} first established the inverse problem of determining time-independent bounded coefficient $q$ when $A=0$ in \eqref{Convection_diffusion_operator}, in a bounded domain, from the full DN map. 
In \cite{Choulli_Stability_potential}, the authors derived logarithmic stability estimates for corresponding problems. Cheng and Yamamoto in \cites{Global_uniquess_cheng_Yamamoto2000,cheng2002identification,cheng_Yamamoto-2004_uniquenss}, determine the time-independent convection term $A$ when $q=0$ in \eqref{Convection_diffusion_operator}, with single boundary measurement in dimension two. For the case of time-independent coefficient, Avdonim and Seidman in \cite{Sergei_Thomas_boundary_contrl_for_q} used boundary control method pioneered by Belishev in \cite{Belishev_Recent_progress_boundary_control} (see also \cite{Kurylev_Book}) to studied unique determination of density coefficient $q$ appearing in \eqref{Convection_diffusion_operator}. In \cite{Bellassoued_full_stably}, Bellassoued and Rassas showed that first-order and zeroth-order time-independent coefficients appearing in a convection-diffusion equation stably recovered from the full DN map.
In \cite{sahoo2019partial}, Sahoo and Vashisth proved recovery of convection and time-dependent density coefficients up to the gauge from Dirichlet data is given on the whole lateral boundary, but Neumann data given only on slightly more than half of the boundary. We also refer to \cite{Bellassoued_partial_stably}, where Bellassoued and Fraj proved the logarithmic stability estimates in determining time-dependent coefficients appearing in the convection-diffusion equation from the partial DN map. Fan and Duan in \cite{fan2021determining} prove unique recovery of time-dependent density coefficient $q$ assuming $A=0$ in \eqref{Convection_diffusion_operator} from partial DN map. Recently, \cite{Anamika-24} establishes unique recovery of time-dependent coefficients $A$ and $q$ up to gauge transform from full Dirichlet and partial Neumann data. Recently, in \cite{mishra2022inverse}, the authors obtained a unique determination of time-dependent coefficients $A$ and $q$ up to gauge transform from full Dirichlet and partial Neumann data in Riemannian geometry setting.

A similar work has been studied by Isakov \cite{Isakov-2007} based on the reflection approach for an inverse problem for Schr\"odinger equation in a bounded domain, where the boundary data are given on an arbitrary part of the boundary whose complement is either a part of hyperplane or sphere. This reflection approach has been used widely in many previous works, for example, Maxwell's systems \cite{Caro-Maxwell-I}, and for the elliptic operators; see, \cites{krupchyk2012inverse, Yang-Paper, Heck_Wang_Optimal_Stability, Anupam-Stability, Bhattacharyya_Optimal_Stability, Bhattacharyya-Kumar, Liu_2024}. In this article, we extend the result obtained by Isakov \cite{Isakov-2007} for the Schr\"odinger equation to the convection-diffusion equation.

\section{Proof of the main theorem}\label{proof of main thm}
\subsection{Preliminary results}\label{Prelim-results}
In this subsection, we shall recall some preliminary results, which we needed to state the existence of Complex Geometric Optics (CGO) solutions to the equation $\mathcal{L}_{A,q}u=0$ and $\Lc^*_{A,q}v=0$ in $(0, T)\times W$, where $W \subset \R^n, n\geq 3$, be bounded open set with smooth boundary.
We shall construct CGO solutions to the equation $\mathcal{L}_{A,q}u=0$ in $(0, T)\times W$ of the form
\begin{equation}
      u(t,x;h)= e^{\frac{\varphi(x)+\mathrm{i}\psi(x)}{h}\eta(t;h)}(a(t,x;h) + r(t,x;h)),
\end{equation}
where $\eta(t;h):= \sin{(h^{2/5}(T-t)^2)}$, $\varphi$ is a limiting Carleman weight for the semiclassical Laplacian (for definition see; \cite{kenig2007calderon}), the real-valued phase function $\psi$ solves the Eikonal equation 
\begin{equation}\label{Eikonal-eqn}
   |\nabla \psi|^{2}=|\nabla \varphi|^{2}, \quad \nabla\varphi\cdot \nabla\psi =0 \text { in } \widetilde{W},
\end{equation}
$a$ is a smooth amplitude, $r$ is a correction term, and $h>0$ is a small semiclassical parameter.

We state the following results for the existence of CGO solutions to the equations $\mathcal{L}_{A,q}u=0$ in $(0,T)\times W$ and $\Lc^*_{A,q}v=0$ in $(0, T)\times W$. We refer readers to see \cite{Anamika-24} for detailed discussion and proofs.
\begin{proposition}[\cite{Anamika-24}]\label{CGO-Soln-prop}
 Let $A\in C^{\infty}(\overline{(0, T)\times W}; \C^{n})$ and $q\in L^{\infty}((0, T)\times W)$. Let $\varphi$ be a limiting Carleman weight, and $\psi$ solves \eqref{Eikonal-eqn} and $\eta$ defined above. Then, for all $h>0$ sufficiently small enough, there exists $h_0>0$ and $r \in L^2(0,T;H^1(W))\cap H^{1}(0,T;H^{-1}(W))$ such that
\begin{align*}
     h\|r\|_{L^2(0,T;H^1_{scl}(W))}\leq \Oc(h^{7/5}),
\end{align*}
and
\begin{align}\label{Form of soln-u}
    u(t,x;h)= e^{\frac{\varphi(x)+\mathrm{i}\psi(x)}{h}\eta(t;h)}(a(t,x;h) + r(t,x;h)),
\end{align}
 is a solution of $\mathcal{L}_{A,q}u=0$, when $h\leq h_0$, and satisfying $u(0,x)=0$, where
 $a(t,x;h)= m(t)e^{\Phi(t,x;h)}$, $m\in C_c^{\infty}(0,T)$ and $\Phi(t,x;h)$ is a solution of 
\begin{align}\label{Operator_Phi_equation}
 \Delta (\varphi(x)+\mathrm{i}\psi(x))+2\nabla(\varphi(x)+\mathrm{i}\psi(x))\cdot\nabla \Phi(t,x;h)+ 2A(t,x)\cdot\nabla(\varphi(x)+\mathrm{i}\psi(x))=0.
 \end{align}
 \end{proposition}
 
\begin{proposition}[\cite{Anamika-24}]\label{Adjoint_CGO-Soln-prop}
 Let $A\in C^{\infty}(\overline{(0, T)\times W}; \C^{n})$ and $q\in L^{\infty}((0, T)\times W)$. Let $\varphi$ be a limiting Carleman weight, and $\psi$ solves \eqref{Eikonal-eqn} and $\eta$ defined above. Then, for all $h>0$ sufficiently small enough, there exists $h_0>0$ and $\widetilde{r} \in L^2(0,T;H^1(W))\cap H^{1}(0,T;H^{-1}(W))$ such that
\begin{align*}
     h\|\widetilde{r}\|_{L^2(0,T;H^1_{scl}(W))}\leq \Oc(h^{7/5}),
\end{align*}
and
\begin{align}\label{Soln-for-adjoint-v}
    v(t,x;h)= e^{\frac{-\varphi(x)+\mathrm{i}\psi(x)}{h}\eta(t;h)}(\widetilde{a}(t,x;h) + \widetilde{r}(t,x;h)),
\end{align}
 is a solution of $\mathcal{L}_{A,q}^{*}v=0$, when $h\leq h_0$, and satisfying $v(T,x)=0$, where 
 $\widetilde{a}(t,x;h)= \widetilde{m}(t)e^{\widetilde{\Phi}(t,x;h)}$, $\widetilde{m}\in C_c^{\infty}(0,T)$ and $\widetilde{\Phi}(t,x;h)$ is a solution of 
 \begin{align}\label{Adjoint-operator_Phi_equation}
\Delta(-\varphi(x)+\mathrm{i}\psi(x))+2\nabla(-\varphi(x)+\mathrm{i}\psi(x))\cdot\nabla \widetilde{\Phi}(t,x;h)- 2\overline{A}(t,x)\cdot\nabla(-\varphi(x)+\mathrm{i}\psi(x))=0.
 \end{align}
\end{proposition}

\subsection{Construction of special solutions}\label{Const_of_special_soln}
In this subsection, we shall construct special solutions $u_{2}$ and $v$ to the equations $\mathcal{L}_{A^{(2)},q_{2}}u_{2}=0$ and $\mathcal{L}_{A^{(1)},q_{1}}^{\ast}v=0$ in $(0,T)\times \Omega$ with 
\begin{equation}\label{Inaccessible-conditions}
    u_{2}|_{(0,T)\times \Gamma_{0}}=0 \mbox { and } v|_{(0,T)\times \Gamma_{0}}=0,
\end{equation}
by utilizing the reflection arguments from \cite{Isakov-2007}. 
In order to achieve the inaccessible Dirichlet conditions in \eqref{Inaccessible-conditions}, we reflect $\Omega$ with respect to the hyperplane $x_{n} = 0$ and denote this reflection by
 \[
 \Omega^{\ast} = \{ (x^{\prime}, -x_{n}):  x =(x^{\prime},x_{n})\in \Omega\}\quad 
\text{ where } \quad x^{\prime} = \left(x_{1}, \cdots , x_{n-1}\right)\in\R^{n-1}.
 \]
Through out this paper, we shall denote $O = \Omega\cup\Omega^*$, $x^{\ast}=(x^{\prime},-x_{n})$ for any $x=(x^{\prime},x_{n})\in \R^{n}$, $f^{\ast}(x)=f(x^{\ast})$ for any function $f$ and $\zeta^*=(\zeta',-\zeta_n)$ for any vector $\zeta=(\zeta',\zeta_n) \in \Cb^n$.

For $k=1,2$ we extend the coefficients $A^{(k)}$ and $q_{k}$ to $(0, T)\times \Omega^{\ast}$ as follows: for the components $A^{(k)}_{j}$, $1\leq j\leq n-1$ and $q_{k}$, we extend them as even functions with respect to the hyperplane $x_{n} = 0$, that is, we define
$$
\begin{array}{rl} \vspace{1ex}
\widetilde{A^{(k)}_{j}}(t, x) := &
\left\{
\begin{array}{ll} \vspace{1ex}
A^{(k)}_{j}(t, x) & \text{ if } (t,x)\in (0,T)\times \Omega \\
A^{(k)}_{j}(t,x^{\ast}) & \text{ if } (t,x)\in (0,T)\times \Omega^{\ast}\\
\end{array}
\right.
, \quad j=1,\cdots,n-1, \\ \vspace{1ex}
\
\widetilde{q_{k}}(t,x):=&
\left\{
\begin{array}{ll} \vspace{1ex}
q_{k}(t, x) &\text{ if } (t,x)\in (0,T)\times \Omega\\
q_{k}(t,x^{\ast}) &\text{ if } (t,x)\in (0,T)\times \Omega^{\ast}\\
\end{array}
\right. .\\
\end{array}
$$
For $A^{(k)}_{n}$, we extend it as an odd function with respect to the hyperplane $x_{n} = 0$; that is, we define
$$
\begin{array}{rl} \vspace{1ex}
\widetilde{A^{(k)}_{n}}(t, x)=&
\left\{
\begin{array}{ll} \vspace{1ex}
A^{(k)}_{n}(t,x) & \text{ if } (t,x)\in (0,T)\times \Omega\\
-A^{(k)}_{n}(t,x^{\ast}) & \text{ if } (t,x)\in (0,T)\times \Omega^{\ast} \\
\end{array}
\right.. \\
\end{array}
$$
Note that $\widetilde{A^{(k)}}\in (C^{\infty}(\overline{(0,T)\times O)})$ and $\widetilde{q_{k}}\in L^{\infty}((0,T)\times O)$, for $k=1,2$.

For $h>0$ is sufficiently small enough, and for any $\xi\in\R^n$ and $\mu^{(1)},\mu^{(2)}\in\R^n$ be such that $|\mu^{(1)}|=|\mu^{(2)}|=1$ and $\mu^{(1)}\cdot\mu^{(2)}=\mu^{(1)}\cdot\xi=\mu^{(2)}\cdot\xi=0$, we choose 
\begin{align*}
\varphi(x) = x \cdot \mu^{(2)} \quad
\psi_{\pm}(x) = x \cdot \bigg[\pm \frac{h^{3/5}\xi}{2} + \sqrt{1-\frac{h^{6/5}|\xi|^2}{4}} \mu^{(1)}\bigg]
\end{align*}
which solves the Eikonal equation \eqref{Eikonal-eqn}.
By Proposition \ref{CGO-Soln-prop}, there exists a CGO solution $\widetilde{u_{2}}\in L^2(0,T;H^1(O))\cap H^{1}(0,T;H^{-1}(O))$ to the equation $\mathcal{L}_{\widetilde{A^{(2)}},\widetilde{q_{2}}}\widetilde{u_{2}}=0$ in $(0,T)\times O$ of the form 
\begin{equation*}
    \widetilde{u_{2}}(t,x,h)= e^{\frac{\varphi(x)+\mathrm{i}\psi_{+}(x)}{h}\eta(t;h)}( m_2(t)e^{\Phi_{2}(t,x;h)}+r_{2}(t,x;h)),
\end{equation*}
where $m_2\in C^{\infty}_c(0,T)$ and $\Phi_{2}\in C^{\infty}((0,T)\times O)$, solves the following transport equation
\begin{equation}\label{Phi_2_equation}
     \Delta(\varphi(x)+\mathrm{i}\psi_{+}(x))+2\left(\nabla\varphi(x)+\mathrm{i}\nabla\psi_{+}(x)\right)\cdot \left(\nabla \Phi_{2}(t,x;h)+\widetilde{A^{(2)}}(t,x)\right) =0, \text { in } (0,T)\times O,
\end{equation} 
and the remainder term $r_{2}\in L^2(0,T;H^1(O))\cap H^{1}(0,T;H^{-1}(O))$, having the decay estimate
\begin{equation}
   h\|r_{2}\|_{L^2(0,T;H^1_{scl}(O))}\leq \Oc(h^{7/5}), \quad \text { as } h\to 0.
\end{equation}
where $h>0$ is a small semiclassical parameter.

\noindent Therefore, we have a CGO solution to the equation $\mathcal{L}_{\widetilde{A^{(2)}},\widetilde{q_{2}}}\widetilde{u_{2}}=0$ in $(0, T)\times O$ of the form 
\begin{equation}\label{CGO-soln-1}
    \widetilde{u_{2}}(t,x,h)= e^{\frac{x\cdot \zeta_{2}}{h}\eta(t;h)}(m_2(t)e^{\Phi_{2}(t,x;h)}+r_{2}(t,x;h)),\quad \text { where } \quad \zeta_2\cdot x :=\varphi(x)+\mathrm{i}\psi_+(x).
\end{equation}
Similarly, by Proposition ~\ref{Adjoint_CGO-Soln-prop}, there exists a CGO solutions $\widetilde{v}\in L^2(0,T;H^1(O))\cap H^{1}(0,T;H^{-1}(O))$ to the equation $\mathcal{L}_{\widetilde{A^{(1)}},\widetilde{q_{1}}}^{*}\widetilde{v}=0$ in $(0,T)\times O$ of the form 
\begin{align*}
    \widetilde{v}(t,x,h)= e^{\frac{-\varphi(x)+\mathrm{i}\psi_{-}(x)}{h}\eta(t;h)}( m_1(t)e^{\Phi_{1}(t,x;h)}+r_{1}(t,x;h)),
\end{align*} 
where $m_1\in C^{\infty}_c(0,T)$ and $\Phi_{1}\in C^{\infty}((0,T)\times O)$, solves the following transport equation
\begin{align}\label{Phi_1_equation}
    \Delta(-\varphi(x)+\mathrm{i}\psi_{-}(x))+2\left(-\nabla\varphi(x)+\mathrm{i}\nabla\psi_{-}(x)\right)\cdot \left(\nabla \Phi_{1}(t,x;h)-\overline{\widetilde{A^{(1)}}}(t,x)\right) =0, \text { in } (0,T)\times O,
\end{align}
and the remainder term $r_{1}\in L^2(0,T;H^1(O))\cap H^{1}(0,T;H^{-1}(O))$, having the decay estimate
\begin{equation}
   h\|r_{1}\|_{L^2(0,T;H^1_{scl}(O))}\leq \Oc(h^{7/5}), \quad \text { as } h\to 0.
\end{equation}

Therefore, we have a CGO solution to the equation $\mathcal{L}_{\widetilde{A^{(1)}},\widetilde{q_{1}}}^{*}\widetilde{v}=0$ in $(0, T)\times O$ of the form  
\begin{equation}\label{CGO-soln-2}
    \widetilde{v}(t,x,h)= e^{\frac{x\cdot \zeta_{1}}{h}\eta(t;h)}(m_1(t)e^{\Phi_{1}(t,x;h)}+r_{1}(t,x;h)), \quad \text { where } \quad \zeta_1\cdot x:=-\varphi(x)+\mathrm{i}\psi_{-}(x).
\end{equation}
From \eqref{Phi_2_equation} and \eqref{Phi_1_equation}, as $h \to 0$, we get the following transport equations
\begin{equation}
\label{eq: Equation for phi_2}
    (\mathrm{i}\mu^{(1)}+\mu^{(2)})\cdot\nabla \Phi_{2}(t,x)+ (\mathrm{i}\mu^{(1)}+\mu^{(2)})\cdot \widetilde{A^{(2)}}(t,x) =0 \quad \text { in } (0,T)\times O,
\end{equation}
and 
\begin{equation}
\label{eq: Equation for phi_1}
    (\mathrm{i}\mu^{(1)}-\mu^{(2)})\cdot\nabla \Phi_{1}(t,x)- (\mathrm{i}\mu^{(1)}-\mu^{(2)})\cdot \overline{\widetilde{A^{(1)}}}(t,x) =0 \quad \text { in } (0,T)\times O.
\end{equation}
Now we can construct special CGO solutions in $(0, T)\times \Omega$, which are vanishing on an inaccessible part of the boundary, namely $(0, T)\times \Gamma_{0}$. To do so, we set
\begin{equation}
    u_{2}(t,x;h):=\widetilde{u_{2}}(t,x;h)-\widetilde{u_{2}}^{\ast}(t,x;h), \quad \text {and} \quad  v(t,x;h):=\widetilde{v}(t,x;h)-\widetilde{v}^{\ast}(t,x;h), 
\end{equation}
where $\widetilde{u}_{2}$ and $\widetilde{v}$ are the CGO solutions of the form \eqref{CGO-soln-1} and \eqref{CGO-soln-2}, respectively. By direct computations, one can easily check that $u_{2}$ and $v$ solve the equations $\mathcal{L}_{A^{(2)},q_{2}}u_{2}=0$ in $(0,T)\times \Omega$ and $\mathcal{L}_{A^{(1)},q_{1}}^{\ast}v=0$ in $(0,T)\times \Omega$, with required inaccessible conditions \eqref{Inaccessible-conditions}. For later purposes, we write down the forms of the special solutions $u_{2}$ and $v$
\begin{align}\label{CGO-u_{2}}
    u_{2}(t,x;h)=e^{\frac{x\cdot \zeta_{2}}{h}\eta(t;h)}(a_{2}(t,x;h)+r_{2}(t,x;h))-e^{\frac{x^{\ast}\cdot \zeta_{2}}{h}\eta(t;h)}(a^{\ast}_{2}(t,x;h)+r^{\ast}_{2}(t,x;h)),
\end{align}
and 
\begin{align}\label{CGO-v}
    v(t,x;h)=e^{\frac{x\cdot \zeta_{1}}{h}\eta(t;h)}(a_{1}(t,x;h)+r_{1}(t,x;h))-e^{\frac{x^{\ast}\cdot \zeta_{1}}{h}\eta(t;h)}(a^{\ast}_{1}(t,x;h)+r^{\ast}_{1}(t,x;h)),
\end{align}
where $a_{j}\in C^{\infty}(\overline{(0,T)\times O})$ and 
$ h\|r_{j}\|_{L^2(0,T;H^1_{scl}(O))}\leq \Oc(h^{7/5}), \quad \text { for } j=1,2.$

\subsection{An integral identity}\label{Integral_identity}
In this subsection, we will use the local DN map to derive an integral identity, which relates the difference of coefficients to the difference of the associated local DN map. 
 Let $u_i$, for $i = 1,2$, be the solutions to the following initial boundary value problem
\begin{equation} \label{eq:4.2}
  \begin{array}{cc}
  \left\{\begin{array}{ r l r}
\mathcal{L}_{A^{(i),q_i}}u_i(t,x) &= 0, & (t,x)\in Q, \\
	 u_i(0,x) &= 0,  & x \in \Omega,   \\
	 u_i(t,x) &= f(t,x), & (t,x) \in \Sigma.
	\end{array}\right.
  \end{array}
 \end{equation}
  Let $u(t,x)=(u_{1}-u_{2})(t,x)$. Then $u$ solves
\begin{equation}\label{eq:4.4}
  \left\{
	\begin{array}{ r l c }
\mathcal{L}_{A^{(1)},q_1}u(t,x) &= 2A(t,x)\cdot\nabla_x u_2(t,x) + q^{\sharp}(t,x)u_2(t,x), \quad &(t,x)\in Q, \\
	 u(0,x) &= 0,  \quad &x \in \Omega,   \\
	 u(t,x) &= 0, \quad &(t,x) \in \Sigma,
	\end{array}
	\right.
 \end{equation}
 where
 \begin{align*}
     &A(t,x) = A^{(1)}(t,x) - A^{(2)}(t,x),\\ 
     &q^{\sharp}_i(t,x) = - \nabla_x\cdot A^{(i)}(t,x) - |A^{(i)}(t,x)|^2 + q_i(t,x),\\ 
     &q^{\sharp}(t,x) = q^{\sharp}_2(t,x) - q^{\sharp}_1(t,x).  
\end{align*}
Let $v(t,x)$ be a solution to the following equation 
\begin{align} \label{eq:52}
   \left\{\begin{array}{rlc}
    \mathcal{L}^*_{A^{(1)},q_1} v(t,x) &= 0, &(t,x)\in Q,  \\
       v(T, x) &=0,  & x \in \O.
\end{array}\right.
\end{align}
Now multiplying both sides in \eqref{eq:4.4} by $\overline{v}$ and then integrating over $Q$, we obtain
\begin{align}
 \int_{Q} (2A\cdot\nabla_x u_2 + q^{\sharp}u_2)\overline{v} 
 \mathrm{d}x \mathrm{d}t &= \int_{Q}  \mathcal{L}_{A^{(1)},q_1}u \overline{v}  \mathrm{d}x \mathrm{d}t \\
 &= \int_{Q} u \overline{ \mathcal{L}^{\ast}_{A^{(1)},q_1} v} \mathrm{d}x \mathrm{d}t - \int_{\Sigma} \partial_{\nu}u \overline{v}\mathrm{d}S_x \mathrm{d}t
 +  \int_{\Omega} u(T,x) \overline{v}(T,x)  \mathrm{d}x,\label{eq:4.6}
\end{align}
where we have used $u|_{\Sigma} = 0, u|_{t=0} = 0$ and $A^{(1)} = A^{(2)}$ on $\Sigma$.  Since $v$ satisfy equation \eqref{eq:52} therefore we get
 \begin{align}
 \label{eq:55}
   \int_{Q} (2A\cdot\nabla_x u_2 + q^{\sharp}u_2)\overline{v} \mathrm{d}x \mathrm{d}t = -\int_{\Sigma} \partial_{\nu}u \overline{v}\mathrm{d}S_x \mathrm{d}t.
 \end{align}
 Now we choose $v$ of the form \eqref{CGO-v}, hence with $\mathrm{supp}(v(t,\cdot)|_{\Sigma})\subset \Gamma$, for each $t\in(0, T)$, therefore \eqref{eq:55} becomes
\begin{align*}
   \int_{Q} (2A\cdot\nabla_x u_2 + q^{\sharp}u_2)\overline{v} \mathrm{d}x \mathrm{d}t = -\int_{\Sigma^{\sharp}} \partial_{\nu}u \overline{v}\mathrm{d}S_x \mathrm{d}t.
 \end{align*}
  From \eqref{eq:Equal_DN_map}, we have $\Lambda_{A^{(1)},q_{1}}^{\Sigma^{\sharp}}(f)= \Lambda_{A^{(2)},q_{2}}^{\Sigma^{\sharp}}(f)$ for $f\in\mathcal{K}_{0}(\Sigma)$, with $\mathrm{supp}(f(t,\cdot)|_{\Sigma})\subset \Gamma$ for each $t\in(0, T)$, imply that $\partial_{\nu}u=0$ on $\Sigma^{\sharp}$.  
Thus, we get the following integral identity
\begin{align}\label{final-integral-identity}
   \int_{Q} (2A\cdot\nabla_x u_2+ q^{\sharp}u_2)\overline{v} \mathrm{d}x \mathrm{d}t =  0,
 \end{align}
where $u_2,v\in L^2(0,T;H^1(\Omega))\cap H^{1}(0,T;H^{-1}(\Omega))$ satisfy \eqref{eq:4.2} and \eqref{eq:52} with inaccessible conditions $ u_{2}|_{(0,T)\times \Gamma_{0}}=0 \mbox { and } v|_{(0,T)\times \Gamma_{0}}=0$ respectively.

\subsection{Recovery of Convection coefficient}\label{Recovery_of_Conv_coeff}
We will insert $u_{2}$ and $v$ into the identity \eqref{final-integral-identity}. To this end, we compute 
\begin{equation}\label{Product of exponential}
\begin{aligned}
&\frac{x\cdot\zeta_{2}}{h}\eta+\frac{x\cdot \overline{\zeta}_{1}}{h}\eta=\mathrm{i}x\cdot\xi \frac{\eta}{h^{2/5}}, 
\quad 
\frac{x\cdot\zeta_{2}}{h}\eta+\frac{x^{\ast}\cdot \overline{\zeta}_{1}}{h}\eta=\mathrm{i}x\cdot \xi_{+}\frac{\eta}{h^{2/5}}+2\mu^{(2)}_{n}x_{n}\frac{\eta}{h},\\
&\frac{x^{\ast}\cdot\zeta_{2}}{h}\eta+\frac{x^{\ast}\cdot \overline{\zeta}_{1}}{h}\eta=
\mathrm{i}x^{\ast}\cdot\xi\frac{\eta}{h^{2/5}}, \quad \frac{x^{\ast}\cdot\zeta_{2}}{h}\eta+\frac{x\cdot \overline{\zeta}_{1}}{h}\eta=\mathrm{i}x\cdot \xi_{-}\frac{\eta}{h^{2/5}}-2\mu^{(2)}_{n}x_{n}\frac{\eta}{h},\ \quad \\
& \text { where } \quad \xi_{\pm}=\left(\xi',\pm\frac{2}{h^{3/5}}\sqrt{1-h^{6/5}\frac{|\xi|^{2}}{4}}\mu^{(1)}_{n}\right).
\end{aligned}
\end{equation} 
To eliminate the undesired terms, we assume $\mu^{(2)}_{n}=0$ and $\mu^{(1)}_{n}\neq 0$, so $|\xi_{+}|$ and $|\xi_{-}|\rightarrow\infty$ as $h\rightarrow 0$.
\noindent Substituting the solutions $u_2$ and $v$ given by \eqref{CGO-u_{2}} and \eqref{CGO-v} into \eqref{final-integral-identity}, we obtain
\begin{equation}\label{ID-1}
\begin{aligned}
  \int_{Q} &2A \cdot \left(\left[ (\zeta_2 \frac{\eta}{h} + \nabla)(a_2+r_2)\right] (\overline{a_1}+\overline{r_1}) e^{\mathrm{i}x\cdot\xi\frac{\eta}{h^{2/5}}}+\left[(\zeta^*_2 \frac{\eta}{h} + \nabla) (a^{\ast}_2+r^{\ast}_2)\right](\overline{a^{\ast}_1}+\overline{r^{\ast}_1}) e^{\mathrm{i}x^{\ast}\cdot \xi\frac{\eta}{h^{2/5}}}\right) \mathrm{d}x\ \mathrm{d}t  \\
   &\quad + \int_{Q} q^{\sharp} (a_2+r_2)(\overline{a_1}+\overline{r_1}) e^{\mathrm{i}x\cdot\xi\frac{\eta}{h^{2/5}}} \mathrm{d}x\ \mathrm{d}t + \int_{Q} q^{\sharp} (a^{\ast}_2+r^{\ast}_2)(\overline{a^{\ast}_1}+\overline{r^{\ast}_1}) e^{\mathrm{i}x^{\ast}\cdot\xi\frac{\eta}{h^{2/5}}} \mathrm{d}x\ \mathrm{d}t\\
   =&\int_{Q} 2A\cdot \left(\left[(\zeta_2 \frac{\eta}{h} + \nabla) (a_2 + r_2)\right](\overline{a^{\ast}_1}+\overline{r^{\ast}_1}) e^{\mathrm{i}x\cdot \xi_{+}\frac{\eta}{h^{2/5}}} 
   +\left[(\zeta^*_2 \frac{\eta}{h} + \nabla) (a^{\ast}_2 + r^{\ast}_2)\right](\overline{a_1}+\overline{r_1}) e^{\mathrm{i}x\cdot \xi_{-}\frac{\eta}{h^{2/5}}}\right) \mathrm{d}x\ \mathrm{d}t \\
   &\quad +\int_{Q} q^{\sharp} (a_2+r_2)(\overline{a^{\ast}_1}+\overline{r^{\ast}_1}) e^{ix\cdot \xi_{+}\frac{\eta}{h^{2/5}}} \mathrm{d}x\ \mathrm{d}t +\int_{Q} q^{\sharp} (a^{\ast}_2+r^{\ast}_2)(\overline{a_1}+\overline{r_1})  e^{\mathrm{i}x\cdot \xi_{-}\frac{\eta}{h^{2/5}}} \mathrm{d}x\ \mathrm{d}t.
\end{aligned}
\end{equation}
Multiplying both sides by $h^{3/5}$, letting $h\rightarrow 0$, and using the Riemann--Lebesgue lemma in the right-hand side of \eqref{ID-1} we get
\begin{align}\label{ID-2}
   &(\mathrm{i}\mu^{(1)}+\mu^{(2)})\cdot \int_{0}^{T}\int_{\Omega} A(t,x) m_1(t)m_2(t)e^{\overline{\Phi_1}+\Phi_2} (T-t)^2e^{\mathrm{i}x\cdot\xi (T-t)^2} \mathrm{d}x\ \mathrm{d}t
   \notag \\& + (\mathrm{i}\mu^{(1)}+\mu^{(2)})^*\cdot \int_{0}^{T}\int_{\Omega} A(t,x) m_1(t)m_2(t)e^{\overline{\Phi^{\ast}_1}+\Phi^{\ast}_2} (T-t)^2e^{\mathrm{i}x^{\ast}\cdot\xi (T-t)^2}\ \mathrm{d}x\ \mathrm{d}t = 0,
\end{align}
for any $\xi,\mu^{(1)},\mu^{(2)}\in\R^n$ such that 
\begin{align}\label{condition on mu and xi}
\mu^{(1)}\cdot\mu^{(2)}=\xi\cdot\mu^{(1)}=\xi\cdot \mu^{(2)}=0, \quad 
|\mu^{(1)}|=|\mu^{(2)}|=1, \quad \mu^{(2)}_{n}=0, \quad \mu^{(1)}_{n}\neq 0.
\end{align}
Making a change of variables, we can rewrite \eqref{ID-2} as
\begin{align*}(\mathrm{i}\mu^{(1)}+\mu^{(2)})\cdot \int_{0}^{T}\int_{O} \widetilde{A}(t,x) (T-t)^2 m_1(t)m_2(t)e^{\overline{\Phi_1}+\Phi_2} e^{\mathrm{i}x\cdot\xi (T-t)^2} \mathrm{d}x\ \mathrm{d}t = 0,
\end{align*}
where $\widetilde{A}$ is the extended vector field in $\overline{(0,T)\times O}$. Since $(T-t)^2$ is non-vanishing for each $t\in (0,T)$ and $ m_1, m_2\in C_c^{\infty}(0,T)$,
varying them leads to
\begin{align*}(\mathrm{i}\mu^{(1)}+\mu^{(2)})\cdot \int_{O} \widetilde{A} e^{\overline{\Phi_1}+\Phi_2} e^{\mathrm{i}x\cdot\xi (T-t)^2} \mathrm{d}x = 0, 
\end{align*}
for each $t\in(0,T)$. 
Without loss of generality, take $\xi(T-t)^2$ as $\xi$. We remove $e^{\overline{\Phi_1}+\Phi_2}$ in the above integral by using a similar idea as \cite{krupchyk2014uniqueness}*{Proposition 3.3}, we obtain
\begin{align}\label{equation_for_A_with_mu}
   (\mathrm{i}\mu^{(1)}+\mu^{(2)})\cdot \int_{O} \widetilde{A} e^{\mathrm{i}x\cdot\xi} \mathrm{d}x = 0, 
\end{align}
for each $t\in (0,T)$ and $\xi,\mu^{(1)},\mu^{(2)}\in\R^n$ with \eqref{condition on mu and xi} holds.
By replacing $\mu^{(2)}$ by $-\mu^{(2)}$, it follows that 
\begin{align}
   \mu\cdot \int_{O}\widetilde{A}  e^{\mathrm{i}x\cdot\xi} \mathrm{d}x =0, 
\end{align}
for each $t\in (0,T)$ and for all $\xi,\mu\in \R^{n}$ such that $\mu\cdot\xi=0$ holds.
Since $A^{(1)}= A^{(2)}$ on $\Sigma$, therefore extending $\widetilde{A}$ by zero on $(0,T)\times (\R^n\setminus O)$, we obtain
\begin{align}\label{equation_for_A}
   \mu\cdot \int_{\R^n}\widetilde{A}  e^{\mathrm{i}x\cdot\xi} \mathrm{d}x = 0.
\end{align}
Let $\xi = (\xi_1,\dots,\xi_n)\in \mathbb{R}^n$, $n\geq 3$, be an arbitrary vector. For an arbitrarily fixed $j,k = 1,\cdots,n$; $j\neq k$, we choose $\mu \in \mathbb{R}^n$ such that $\mu_k=-\xi_j$, $\mu_j = \xi_k$ and $\mu_{l}=0$ whenever $l\neq j,k$.
Then $\mu \cdot \xi =0$, and therefore the above identity becomes
\begin{align}
    \xi_j\widehat{\widetilde{A_{k}}}(t,\xi) - \xi_k\widehat{\widetilde{A_{j}}}(t,\xi) = 0, \quad \text { in } \R^n, \quad \text { for each } t\in (0,T).
\end{align}
 Hence, in particular, $\partial_{x_j}\widetilde{A_{k}}(t,x)-\partial_{x_k}\widetilde{A_{j}}(t,x)=0  \text{ in } \overline{(0,T)\times 
 O}$, for $j\neq k$ and thus, $d A=0$, in $(0,T)\times \Omega$. Therefore, using simply-connectedness of $\Omega$, there exists $\Psi\in C^{\infty}(\overline{Q})$ with $\Psi|_{\Sigma} = 0$ such that $A^{(2)}(t,x)-A^{(1)}(t,x)=\nabla_{x}\Psi(t,x)$ in $Q$. Moreover, we get $\p_{\nu}\Psi|_{\Sigma} = 0$ since $A^{(1)}=A^{(2)}$ on $\Sigma$.

\subsection{Recovery of density coefficient}\label{Recov_of_dens_coeff}
In this subsection, we establish the uniqueness of the density coefficient up to gauge invariance. From Subsection ~\ref{Recovery_of_Conv_coeff}, we know that $A^{(2)}(t,x)-A^{(1)}(t,x)=\nabla_{x}\Psi(t,x)$ for some $\Psi\in C^{\infty}(\overline{Q})$ with $\Psi|_{\Sigma} = \partial_{\nu}\Psi|_{\Sigma} = 0$. To proceed, we replace the pair $(A^{(1)},q_1)$ by $(A^{(3)},q_3)$ where $A^{(3)} = A^{(1)} + \nabla_x\Psi$ and $q_3 = q_1 + \partial_t\Psi$. Consequently, we have $\Lambda_{A^{(3)},q_3}^{\Sigma^{\sharp}} = \Lambda_{A^{(2)},q_2}^{\Sigma^{\sharp}}$ (we refer to Remark \ref{Gauge-Invariance}). Using this fact along with $A^{(2)} = A^{(3)} $, we derive the following integral identity
\begin{align*}
     \int_{Q} (q_2-q_3)(t,x)u_2(t,x)\overline{v(t,x)}\mathrm{d}x\mathrm{d}t = 0.
\end{align*}
Substituting the solutions $u_2$ and $v$ as given by \eqref{CGO-u_{2}} and \eqref{CGO-v} into the above equation, and then taking $h \rightarrow 0$ while applying the Riemann--Lebesgue lemma, we obtain
\begin{align*}
     \int_{0}^{T}\int_{\Omega} q(t,x)m_1(t)m_2(t) \left( e^{\overline{\Phi_1(t,x)}+\Phi_2(t,x)}e^{\mathrm{i}x\cdot\xi(T-t)^2} +  e^{\overline{\Phi^*_1(t,x)}+\Phi^*_2(t,x)}e^{\mathrm{i}x^{\ast}\cdot\xi (T-t)^2} \right)\mathrm{d}x\mathrm{d}t = 0,
\end{align*}
 where $q(t,x) = (q_2-q_3)(t,x) $. Applying the change of variables, we get
\begin{align*}
     \int_{0}^{T}\int_{O} \widetilde{q}(t,x) m_1(t)m_2(t) e^{\overline{\Phi_1(t,x)}+\Phi_2(t,x)}e^{\mathrm{i}x\cdot\xi(T-t)^2}\mathrm{d}x\mathrm{d}t=0,
\end{align*}
where $\widetilde{q}$ is the extended function in $(0,T)\times O$.
We can replace $e^{\Phi_2}$ with $e^{\Phi_2}g$ in the expression for $u_2$, provided that $g\in C^{\infty}(\overline{(0, T)\times O})$ and satisfies $(\mathrm{i}\mu^{(1)}+\mu^{(2)})\cdot\nabla g = 0$ in $(0, T)\times O$. Thus, the above equation can be replaced by
\begin{align}
\label{eq: Equation for g(t,x)}
\int_{0}^{T}\int_{O}\widetilde{q}(t,x) m_1(t)m_2(t)g(t,x) e^{\overline{\Phi_1(t,x)}+\Phi_2(t,x)}e^{\mathrm{i}x\cdot\xi(T-t)^2}\mathrm{d}x\mathrm{d}t  = 0.
\end{align}
By combining Equations \eqref{eq: Equation for phi_2} and \eqref{eq: Equation for phi_1} with the equality $A^{(2)} = A^{(3)}$, we get
\begin{align*}
  (\mathrm{i}\mu^{(1)}+\mu^{(2)})\cdot\nabla (\overline{\Phi_{1}(t,x)}+ \Phi_{2}(t,x)) =0 \quad \text { in } (0,T)\times O.
\end{align*}
Letting, $\displaystyle g(t,x) = e^{-\left(\overline{\Phi_{1}(t,x)} + \Phi_{2}(t,x)\right)}$ and substituting this into \eqref{eq: Equation for g(t,x)}, we obtain
\begin{align*}
\int_{0}^{T}\int_{O}\widetilde{q}(t,x) m_1(t)m_2(t)e^{\mathrm{i}x\cdot\xi(T-t)^2}\mathrm{d}x\mathrm{d}t  = 0.
\end{align*}
Since $(T-t)^2$ is non-vanishing for $t\in (0,T)$ and varying $m_1, m_2 \in C_c^{\infty}(0,T)$ and denoting $\xi(T-t)^2$ as $\xi$, we get
\begin{align}\label{fourier_of_q}
\int_{O}\widetilde{q}(t,x) e^{\mathrm{i}x\cdot\xi}\ \mathrm{d}x  = 0,
\end{align}
for each $t\in (0,T)$ and for all $\xi,\mu^{(1)},\mu^{(2)}\in\R^n$ with \eqref{condition on mu and xi} holds.
For any $\xi\in\mathbb{R}^n$ we can find $\mu^{(1)}$ and $\mu^{(2)}$ such that \eqref{condition on mu and xi} holds, then from \eqref{fourier_of_q} we have $\widehat{\widetilde{q}}(t,\xi) = 0$ for each $t\in(0, T)$. 
Therefore, we conclude that $q_2(t,x)-q_1(t,x) = \partial_t\Psi(t,x)$ in $Q$. This completes the proof of Theorem ~\ref{Theorem_hyperplane}.

\section*{acknowledgments}
The authors thank Sombuddha Bhattacharyya and Rohit Kumar Mishra for several fruitful discussions throughout the project. PK would like to thank the Indian Institute of Technology, Gandhinagar, for the hospitality during his visit, where a part of this work was completed. AP  is supported by UGC, Government of India, with a research fellowship.

\bibliographystyle{plain}
\bibliography{bibfile.bib}
\end{document}